# A STRONG INVARIANCE PRINCIPLE FOR ASSOCIATED RANDOM FIELDS[??]

By Raluca M. Balan

*University of Ottawa*

In this paper we generalize Yu's [*Ann. Probab.* **24** (1996) 2079–2097] strong invariance principle for associated sequences to the multi-parameter case, under the assumption that the covariance coefficient $u(n)$ decays exponentially as $n \to \infty$. The main tools that we use are the following: the Berkes and Morrow [*Z. Wahrsch. Verw. Gebiete* **57** (1981) 15–37] multi-parameter blocking technique, the Csörgő and Révész [*Z. Wahrsch. Verw. Gebiete* **31** (1975) 255–260] quantile transform method and the Bulinski [*Theory Probab. Appl.* **40** (1995) 136–144] rate of convergence in the CLT.

**1. Introduction.** Among various concepts introduced to measure the dependence between random variables, association deserves a special place because of its numerous applications and its relatively easy mathematical manipulation. A finite collection $(X_1, \ldots, X_m)$ of random variables is said to be *associated* (or satisfies the *FKG inequalities*) if for any coordinatewise nondecreasing functions $f, g$ on $\mathbf{R}^m$, $\mathrm{cov}(f(X_1, \ldots, X_m), g(X_1, \ldots, X_m)) \geq 0$, whenever the covariance is defined. An infinite collection of random variables is associated if every finite sub-collection is associated. This concept was formally introduced by Esary, Proschan and Walkup (1967), who also deduced some of its most important properties.

In the past few decades, a lot of effort has been dedicated to prove limit theorems for random fields $(X_j)_{j \in \mathbf{Z}_+^d}$ of associated random variables. In the case $d = 1$, this culminated with the strong invariance principle of Yu (1996), from which one can easily deduce all the other major limit theorems, like the weak invariance principle and the functional law of the iterated logarithm

Received November 2003; revised May 2004.

[1]Supported in part by a research grant from the Natural Sciences and Engineering Research Council of Canada.

(FLIL). The present paper was motivated by the need for a similar result in the case $d \geq 2$, which arises in the context of higher-dimensional models, like the percolation model of Cox and Grimmett (1984).

The first asymptotic result for zero-mean associated random fields was the central limit theorem (CLT) proved by Newman (1980) for the (strongly) stationary case. This result says that if the *finite susceptibility* assumption holds, that is, $\sigma^2 := \sum_{i \in \mathbf{Z}^d} \rho(i) < \infty$, where $\rho(j-k) := \mathrm{cov}(X_j, X_k)$, then

$$n^{-d/2} S_n \xrightarrow{d} N(0, \sigma^2), \tag{1}$$

where $S_n := \sum_{j_1 \leq n} \cdots \sum_{j_d \leq n} X_j$. This was generalized by Cox and Grimmett (1984) to the nonstationary case, under the assumption that $u(n) \to 0$ as $n \to \infty$, where

$$u(n) := \sup_{j \in \mathbf{Z}_+^d} \sum_{k : \|j-k\| \geq n} \mathrm{cov}(X_j, X_k) \tag{2}$$

and $\|i\| := \max_{s=1,\ldots,d} |i_s|$.

The weak invariance principle for (strongly) stationary associated random fields was proved by Newman and Wright (1981, 1982) in the case $d=1$ and $d=2$, under the same finite susceptibility assumption. They also conjectured that the same principle holds for $d > 2$. A partial solution to this problem was given by Burton and Kim (1988) in the stationary case, and by Kim (1996) in the nonstationary case, under the *finite $r$-susceptibility* assumption:

$$E|S_N|^{2+r} \leq C[N]^{1+r/2},$$

where $S_N := \sum_{j \leq N} X_j$ and $[N] := \prod_{s=1}^d N_s$ for $N = (N_1, \ldots, N_d) \in \mathbf{Z}_+^d$. (If $i, j \in \mathbf{Z}_+^d$, we use the notation $i \leq j$ if $i_s \leq j_s$, $\forall s = 1, \ldots, d$ and $i < j$ if $i_s < j_s$, $\forall s = 1, \ldots, d$.)

The conjecture was fully solved by Bulinski and Keane (1996), who proved that for a zero-mean (weakly) stationary associated random field $(X_j)_{j \in \mathbf{Z}_+^d}$ with uniformly bounded moments of order $s > 2$ and a power decay rate for the covariance coefficient $u(n)$, we have

$$W_n(\cdot) \xrightarrow{d} W(\cdot) \qquad \text{in } D([0,1]^d), \tag{3}$$

where $W_n(t) := n^{-d/2} \sum_{j_1 \leq nt_1} \cdots \sum_{j_d \leq nt_d} X_j$ and $W = (W(t))_{t \in [0,1]^d}$ is a $d$-parameter Wiener process with variance $\sigma^2$. We note, in passing, that for $d = 1$, generalizations to the nonstationary case and to the vector-valued case were given by Birkel (1988a) and by Burton, Dabrowski and Dehling (1988), respectively.

The FLIL for associated sequences was obtained by Dabrowski (1985), under the finite $r$-susceptibility assumption with $r = 1$ and a condition which requires that $E(S_n^2)/n$ converges to $\sigma^2$ with a power decay rate.

The strong invariance principle proved by Yu (1996) strengthened and unified all of these results in the case $d = 1$ and implied other asymptotic fluctuation results, like the Chung's type of FLIL for the maxima of partial sums; see Theorems A–E of Philipp and Stout (1975). More precisely, Yu showed that if $(X_j)_{j \in \mathbf{Z}_+}$ is a sequence of associated random variables such that the moments of order $s > 2$ are uniformly bounded, the variances are bounded below away from 0 and the covariance coefficient $u(n)$ decays exponentially as $n \to \infty$, then it is possible to redefine the original sequence on a richer probability space together with a standard Wiener process $W = (W(t))_{t \in [0, \infty)}$ such that, for some $\varepsilon > 0$

$$S_n - W(\sigma_n^2) = O(n^{1/2-\varepsilon}) \qquad \text{a.s.,}$$

where $\sigma_n^2 := E(S_n^2)$. As far as we know, there are no generalizations of this principle to the case $d \geq 2$. The purpose of the present paper is to fill this gap and to provide a powerful approximation tool that can be used in higher dimensions.

Unlike the case $d = 1$, the strong invariance principle for associated random fields in higher dimensions holds only for points $N \in \mathbf{Z}_+^d$ which are not "too close" to the coordinate axes. This is not at all surprising and a similar fact happens for mixing random fields; see Theorem 1 of Berkes and Morrow (1981). The reason for this phenomenon is the irregular behavior of $E(S_N^2)$ close to the coordinate planes.

We proceed now to introduce the notation that will be used throughout this paper.

Let $(X_j)_{j \in \mathbf{Z}_+^d}$ be a weakly stationary associated random field with zero mean and $\rho(j - k) := E(X_j X_k)$, $\forall j, k \in \mathbf{Z}_+^d$. Let $u(n)$ be the covariance coefficient defined by (2). Because of stationarity, we have $u(n) = \sum_{i \in \mathbf{Z}^d : \|i\| \geq n} \rho(i)$ for every $n \geq 0$. We will suppose that $\rho(0) > 0$ and $\sigma^2 := u(0) = \sum_{i \in \mathbf{Z}^d} \rho(i) < \infty$.

For any finite subset $V \subseteq \mathbf{Z}_+^d$, we let $|V|$ be the cardinality of $V$, $S(V) := \sum_{j \in V} X_j, \sigma^2(V) := E[S^2(V)]$ and $F_V(x) := P(S(V)/\sigma(V) \leq x), x \in \mathbf{R}$. Note that for any finite subset $V \subseteq \mathbf{Z}_+^d$,

(4) $$\rho(0) \leq \frac{\sigma^2(V)}{|V|} \leq \sigma^2.$$

Most of the time we will work with "rectangles" $V \subseteq \mathbf{Z}_+^d$ of the form $V := (a, b] = \prod_{s=1}^d (a_s, b_s]$ with $a_s, b_s \in \mathbf{Z}_+ \cup \{0\}, a_s \leq b_s$; note that $|V| = [b - a]$. We denote with $\mathcal{A}$ the class of all the subsets $V$ of this form.

We will use the following conditions:

(C1) $\sup_{j \in \mathbf{Z}_+^d} E|X_j|^{2+r+\delta} < \infty$ for some $r, \delta > 0$.
(C2) $u(n) = O(e^{-\lambda n})$ for some $\lambda > 0$.

(C2′) $u(n) = O(n^{-\nu})$ for some $\nu > 0$.

We recall that a *d-parameter Wiener process* $W = \{W_t; t \in [0,\infty)^d\}$ with variance $\sigma^2$ is a Gaussian process with independent increments such that $W(R)$ has a $N(0, \sigma^2|R|)$-distribution for any rectangle $R$ ($|R|$ denotes the volume of $R$). Using the same notation as Berkes and Morrow (1981), we put

$$G_\tau := \bigcap_{s=1}^d \left\{ j \in \mathbf{Z}_+^d : j_s \geq \prod_{s' \neq s} j_{s'}^\tau \right\}, \qquad \tau \in (0,1).$$

Here is the main result of this paper.

THEOREM 1.1.   *Let $d \geq 2, \tau \in (0,1)$ and $(X_j)_{j \in \mathbf{Z}_+^d}$ be a weakly stationary associated random field with zero mean and $\rho(j - k) := E(X_j X_k)$ for any $j, k \in \mathbf{Z}_+^d$. Suppose that $\rho(0) > 0$ and $\sigma^2 := \sum_{i \in \mathbf{Z}^d} \rho(i) < \infty$.*

*If* (C1) *and* (C2) *hold, then without changing its distribution, we can redefine the random field $(X_j)_{j \in \mathbf{Z}_+^d}$ on a richer probability space together with a d-parameter Wiener process $\{W_t; t \in [0,\infty)^d\}$ with variance $\sigma^2$ such that*

$$S_N - W_N = O([N]^{1/2 - \varepsilon}) \qquad a.s.$$

*for $N \in G_\tau$. Here $\varepsilon$ is a positive constant depending on the field $(X_j)_{j \in \mathbf{Z}_+^d}$.*

From the previous theorem one can easily deduce the following CLT:

$$[N]^{-1/2} S_N \xrightarrow{d} N(0, \sigma^2)$$

when $[N] \to \infty$ and $N \in G_\tau$ for some $\tau \in (0,1)$; this is more general than (1) which was obtained only for $N = (n, \ldots, n) \in \mathbf{Z}_+^d$. The nonfunctional version of LIL obtained by Wichura (1973) for any multi-parameter process with independent increments (in particular, for the Wiener process) allows us to conclude that

$$\limsup_{[N] \to \infty, N \in G_\tau} (2[N] \log \log [N])^{-1/2} S_N = \sigma \qquad a.s.$$

We proceed now to the proof of Theorem 1.1. This is divided into several steps which are explained in Section 2. The remaining sections contain the developments that are needed to perform each step. To ease the exposition, we placed in the Appendix the proofs of some preliminary lemmas.

**2. Description of the method.** In this section we will indicate what are the main ingredients that are needed for the proof of Theorem 1.1. More precisely, by blending the multi-parameter blocking technique of Berkes and Morrow (1981) with the quantile transform technique of Csörgő and Révész (1975), we will be able to generalize to the multi-parameter case the method introduced by Yu (1996).

Throughout our work we will use the letter $C$ to denote a generic positive constant, independent of $k$.

Let $\alpha > \beta > 1$ be integers to be chosen later and $n_0 := 0$. For $l \in \mathbf{Z}_+$, let

$$n_l := \sum_{i=1}^{l}(i^\alpha + i^\beta) \sim \frac{1}{\alpha+1}l^{\alpha+1}.$$

For each $k := (k_1, \ldots, k_d) \in \mathbf{Z}_+^d$, we put $N_k := (n_{k_1}, \ldots, n_{k_d})$. For all $k \in \mathbf{Z}_+^d$, we have $[N_k] \sim (\alpha+1)^{-d}[k]^{\alpha+1}$.

Let $B_k := (N_{k-1}, N_k] = \prod_{s=1}^{d}(n_{k_s-1}, n_{k_s}]$. Note that $|B_k| = \prod_{s=1}^{d}(k_s^\alpha + k_s^\beta) \leq 2^d[k]^\alpha$. We define the "big" blocks $H_k$ and the "small" blocks $I_k$ by

$$H_k := \prod_{s=1}^{d}(n_{k_s-1}, n_{k_s-1} + k_s^\alpha], \qquad I_k := B_k \setminus H_k.$$

Note that $|H_k| = [k]^\alpha$ and $(2^d-1)[k]^\beta \leq |I_k| \leq (2^d-1)[k]^\alpha$. We denote $u_k := S(H_k), \lambda_k^2 := \sigma^2(H_k)$ and $v_k := S(I_k), \tau_k^2 := \sigma^2(I_k)$. By (4),

(5) $\qquad C[k]^\alpha \leq \lambda_k^2 \leq C[k]^\alpha, \qquad C[k]^\beta \leq \tau_k^2 \leq C[k]^\alpha.$

The sums over the big blocks will be used to generate a Gaussian approximating sequence $(\eta_k)_k$ which will in turn be approximated by a Wiener process. In order to do this, we will need an upper bound for the covariance of the sums over two big blocks in terms of the distance between these blocks. The small blocks are introduced simply to give some space between the big blocks, that is, to ensure that the distance between any two big blocks is nonzero.

If the distribution function $\tilde{F}_k$ of $u_k/\lambda_k$ is continuous, then one could use directly the quantile transform method of Csörgő and Révész (1975) to approximate the variable $u_k/\lambda_k$ by a $N(0,1)$-random variable. In general, this assumption may not be satisfied, and, therefore one, needs to employ a "smoothing" technique [see Yu (1996)]. Without changing its distribution, we redefine the random field $(u_k)_{k\in\mathbf{Z}_+^d}$ on a rich enough probability space together with a random field $(w_k)_{k\in\mathbf{Z}_+^d}$ of independent random variables such that $w_k$ is $N(0, \tau_k^2)$-distributed and $(u_k)_k$ and $(w_k)_k$ are independent. Let

$$\xi_k := (u_k + w_k)/(\lambda_k^2 + \tau_k^2)^{1/2}, \qquad k \in \mathbf{Z}_+^d,$$

and $F_k$ be the distribution function of $\xi_k$. By the CLT for associated random fields, $\tilde{F}_k(x) \to \Phi(x)$ as $k \to \infty$ and, consequently, $F_k(x) \to \Phi(x)$ as $k \to \infty$, where $\Phi(x)$ denotes the $N(0,1)$ distribution function. Therefore, it is reasonable to consider the following $N(0,1)$-random variable

$$\eta_k := \Phi^{-1}(F_k(\xi_k))$$

as an approximation for $\xi_k$. Let $e_k := \sqrt{\lambda_k^2 + \tau_k^2}(\xi_k - \eta_k)$.

In what follows we will adapt the method introduced by Berkes and Morrow (1981) for mixing random fields to suit the special needs of an associated random field.

Following page 25 of Berkes and Morrow (1981), we let $\tau \in (0,1)$ be arbitrary, $\rho := \tau/8$, $L$ be the set of all indices $i$ corresponding to the "good" blocks $B_i \subseteq G_\rho$, and $H$ be the set of all points in $\mathbf{Z}_+^d$ which fall in one of the good blocks. The good blocks collect those points $n \in \mathbf{Z}_+^d$ which are not too close to the coordinate axes. From the technical point of view, their indices $i$ satisfy the relationship $i_s \geq C[i]^{\rho/2}$, which is crucial in the proof of Lemma 3.9.

To each point $N \in H$, we associate the points $N^{(1)}, \ldots, N^{(d)}$ which can be thought as the intersections of the hyperplanes $n_s = N_s, s = 1, \ldots, d$ with the "boundary" of the domain $H$; their precise definition is $N_{s'}^{(s)} = N_{s'}$, $\forall s' \neq s$ and

$$N_s^{(s)} := \min_{n \in H \,:\, n_{s'}=N_{s'}, s' \neq s} n_s.$$

Unlike the above-mentioned authors, we raise a small technical point by noting that $H$ may not be a nice "L-shaped" region. This is why we consider the rectangles $R_k := (M_k, N_k] \subseteq H$, where $M_k := ((N_k^{(1)})_1, \ldots, (N_k^{(d)})_d)$. We note that $L_k := \{i : B_i \subseteq R_k\} \subseteq L \cap \{i \leq k\}$.

If $V$ is a rectangle in $\mathbf{Z}_+^d$ and $\tilde{V}$ is the rectangle in $\mathbf{R}_+^d$ which corresponds to $V$, then we make an abuse of notation by writing $W(V)$ instead of $W(\tilde{V})$.

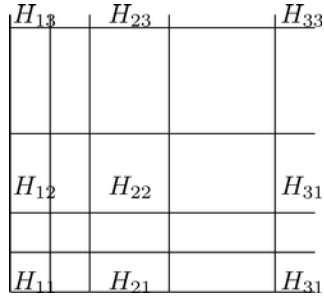

Fig. 1.

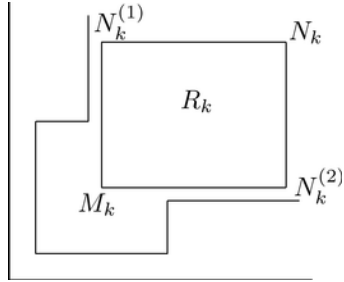

Fig. 2.

This convention will be used throughout this work and will occasionally apply to finite unions of rectangles as well. We write

$$S_N = (S_N - S_{N_k}) + S(R_k) + S((0, N_k] \setminus R_k),$$
$$W_N = (W_N - W_{N_k}) + W(R_k) + W((0, N_k] \setminus R_k)$$

and we use the following decomposition of $S(R_k)$, based on the definitions of $\xi_i$ and $e_i$ and the fact that $S(B_i) = u_i + v_i$:

$$S(R_k) = \sum_{i \in L_k} e_i + \sum_{i \in L_k} \sqrt{|B_i|} \left( \sqrt{\frac{\lambda_i^2 + \tau_i^2}{|B_i|}} - \sigma \right) \eta_i$$
(6)
$$+ \sum_{i \in L_k} \sigma \sqrt{|B_i|} \, \eta_i - \sum_{i \in L_k} w_i + \sum_{i \in L_k} v_i.$$

In Section 3 we will show that all the sums in the above decomposition, except the third one, can be made sufficiently small. The third sum will be treated separately in Section 4 and will be approximated by $W(R_k) = \sum_{i \in L_k} W(B_i)$, via a very powerful approximation result [Theorem 5 of Berkes and Philipp (1979)] and a carefully chosen procedure for counting the indices in $L$. Finally, in Section 5 we will show that the terms $S((0, N_k] \setminus R_k), W((0, N_k] \setminus R_k)$ can be made sufficiently small if $N_k \in G_\tau$, and the differences $S_N - S_{N_k}, W_N - W_{N_k}$ are small if $N \in G_\rho$. This will conclude the proof of Theorem 1.1.

**3. The "good" blocks.** In this section we will show that all the sums in the decomposition (6) of $S_{R_k}$, except the third one, can be made sufficiently small.

In order to treat the first sum of this decomposition, we need to evaluate the precision of the approximation of $\xi_k$ by $\eta_k$. This will be given by the speed of convergence in the CLT. In this paper we decided to use the result obtained by Bulinski (1995), under the assumption that the covariance coefficient $u(n)$ decays exponentially as $n \to \infty$. Under this assumption, this is the sharpest speed of convergence in the CLT when $d = 1, s = 3$ [see Birkel (1988b)]. We note in passing that in the case $d = 1$, a different speed of convergence in the CLT was developed and used by Yu (1996) for associated sequences with a power decay rate of the covariance coefficient; however, the exponential decay rate of $u(n)$ was eventually needed for the strong invariance principle. The problem of whether or not the strong invariance principle continues to hold for associated random fields with a power decay rate of covariances is still open even in the case $d = 1$, and we do not attempt to tackle it here.

LEMMA 3.1 [Theorems 1 and 2 of Bulinski (1995)]. *Suppose that* (C1) *and* (C2) *hold and let* $s := 2 + r + \delta$. *Then for any finite subset* $V \subseteq \mathbf{Z}_+^d$,

$$\sup_{x \in \mathbf{R}} |F_V(x) - \Phi(x)|$$

$$\leq \begin{cases} C|V| \cdot (\sigma^2(V))^{-s/2} \cdot (\log(|V| + 1))^{d(s-1)}, & \text{if } s \leq 3, \\ C|V| \cdot (\sigma^2(V))^{-3/2} \cdot (\log(|V| + 1))^{d}, & \text{if } s > 3. \end{cases}$$

The next result is a generalization of Lemma 3.2 of Yu (1996) to the case $d \geq 2$, in the case of an exponential decay rate of $u(n)$. Its proof is routine and is given in the Appendix.

LEMMA 3.2. *If* (C1) *and* (C2) *hold and* $2r_0 r/(2+r) < \alpha/\beta < 2(1+r)/(2+r)$ *with* $r_0 := \max\{1, (r+\delta)^{-1}\}$, *then*

$$\sup_{x \in \mathbf{R}} |F_k(x) - \Phi(x)| \leq C[k]^{-r\beta/(2+r)} \quad \text{and} \quad \sup_{x \in \mathbf{R}} |f_k(x) - f(x)| \leq C,$$

*where $f_k(x)$ is the density function of $\xi_k$ and $f(x)$ is the $N(0,1)$ density function.*

Using Lemma 3.2 and an argument that was introduced in the proof of Lemma 3 of Csörgő and Révész (1975), we get the precision of the approximation of $\xi_k$ by $\eta_k$.

LEMMA 3.3. *Under* (C1) *and* (C2), *we have*

$$|\Phi^{-1}(F_k(x)) - x| \leq C[k]^{-\{r\beta/(2+r) - K^2/2\}},$$

*provided that* $|x| \leq K\sqrt{\log[k]}$, *where* $0 < K < \sqrt{2r\beta/(2+r)}$.

Next we give the precision of the approximation of $\xi_k$ by $\eta_k$ in terms of the $L^2$-distance. For this we will need the following lemma which gives an upper bound for the moments of order $2 + r$, generalizing an older result of Birkel (1988c) in the case $d = 1$. In particular, this lemma shows that $(X_j)_{j \in \mathbf{Z}_+^d}$ has finite $r$-susceptibility (as defined in the Introduction).

LEMMA 3.4 [Corollary 1 of Bulinski (1993)]. *Suppose that* (C1) *and* (C2$'$) *hold with* $\nu \geq d\nu_0$, *where* $\nu_0 := r(2+r+\delta)/(2\delta) < (d-2)^{-1}$ *if* $d \geq 3$. *Then for any* $V \in \mathcal{A}$,
$$E|S(V)|^{2+r} \leq C|V|^{1+r/2}.$$

Using (5), Lemmas 3.3 and 3.4, and employing the same technique that was used in the proof of Lemma 3.10 of Yu (1996), we get the following result.

LEMMA 3.5. *Under* (C1) *and* (C2), *we have*
$$E[e_k^2] \leq C[k]^{\alpha - \varepsilon_0} \qquad \forall k \in \mathbf{Z}_+^d,$$
*where* $\varepsilon_0 := 2r^2\beta/\{(2+r)(4+3r)\}$.

The next result will show us that the first sum in the decomposition (6) of $S(R_k)$ is small.

LEMMA 3.6. *Suppose that* (C1) *and* (C2) *hold and* $\beta > (1+2/r)(3+4/r)$. *Then there exists* $\varepsilon_1 > 0$ *such that for every* $k \in \mathbf{Z}_+^d$ *with* $L_k \neq \varnothing$,
$$\sum_{i \in L_k} |e_i| \leq C[N_k]^{1/2 - \varepsilon_1} \qquad a.s.$$

PROOF. Let $q > 0$ be such that $\alpha - \varepsilon_0 + 1 < 2q < \alpha - 1$ (this is possible since $\varepsilon_0 > 2$ by our choice of $\beta$). By Chebyshev's inequality and Lemma 3.5, we have
$$P(|e_i| \geq [i]^q) \leq [i]^{-\{2q - (\alpha - \varepsilon_0)\}} \qquad \forall i \in \mathbf{Z}_+^d.$$
By the Borel–Cantelli lemma, it follows that $|e_i| \leq C[i]^q$, $\forall i \in \mathbf{Z}_+^d$ a.s. and, hence, $\sum_{i \in L_k} |e_i| \leq C \sum_{i \in L_k} [i]^q \leq C[k]^{q+1} \leq C[k]^{(\alpha+1)/2 - \varepsilon_1'} \leq C[N_k]^{1/2 - \varepsilon_1}$ a.s., where $0 < \varepsilon_1' < (\alpha - 1)/2 - q$ and $\varepsilon_1 := \varepsilon_1'/(\alpha+1)$. □

The proof of the following lemma is given in the Appendix.

LEMMA 3.7. *If* (C2$'$) *holds with* $d < \nu < 2d$, *then*

(7) $$\sigma^2 - \frac{\sigma^2(V)}{|V|} = O(|V|^{-\delta_0}),$$

*where* $V$ *is a finite union of rectangles in* $\mathcal{A}$ *and* $\delta_0 := \nu/d - 1$.

REMARK. Relationship (7) is exactly Dabrowski's (1985) condition for the FLIL for associated sequences.

LEMMA 3.8. *Suppose that* (C2$'$) *hold with* $d < \nu < 2d$ *and* $\beta > 3/\delta_0$, *where* $\delta_0 := \nu/d - 1$. *Then for every* $k \in \mathbf{Z}_+^d$ *with* $L_k \neq \varnothing$,

$$\sum_{i \in L_k} \sqrt{|B_i|}\left(\sigma - \sqrt{\frac{\lambda_i^2 + \tau_i^2}{|B_i|}}\right)|\eta_i| \leq C[N_k]^{1/2-\alpha_0} \qquad a.s.,$$

*where* $\alpha_0 := 1/\{2(\alpha+1)\}$.

PROOF. Note that $a_i := \sigma - \sqrt{(\lambda_i^2 + \tau_i^2)/|B_i|} > 0$, by (4) and the association property. Using (7), we have

$$a_i^2 \leq \sigma^2 - \frac{\lambda_i^2 + \tau_i^2}{|B_i|} = \frac{|H_i|}{|B_i|}\left(\sigma^2 - \frac{\lambda_i^2}{|H_i|}\right) + \frac{|I_i|}{|B_i|}\left(\sigma^2 - \frac{\tau_i^2}{|I_i|}\right)$$
$$\leq C(|H_i|^{-\delta_0} + |I_i|^{-\delta_0}) \leq C[i]^{-\beta\delta_0}$$

and, hence, by Chebyshev's inequality,

$$P(\sqrt{|B_i|}a_i\eta_i \geq [i]^{\alpha/2-1}) \leq [i]^{-(\alpha-2)}|B_i|a_i^2 \leq C[i]^{-(\beta\delta_0-2)}.$$

By the Borel–Cantelli lemma, it follows that $\sqrt{|B_i|}a_i\eta_i \leq C[i]^{\alpha/2-1}$, $\forall i \in \mathbf{Z}_+^d$ a.s. and, hence, $\sum_{i \in L_k} \sqrt{|B_i|}a_i\eta_i \leq C\sum_{i \in L_k}[i]^{\alpha/2-1} \leq C[k]^{\alpha/2} \leq C[N_k]^{1/2-\alpha_0}$ a.s. since $[k] \sim (\alpha+1)^{d/(\alpha+1)}[N_k]^{1/(\alpha+1)}$. □

The final result of this section shows that the last two sums in the decomposition (6) of $S(R_k)$ are small.

LEMMA 3.9. *If* $\alpha - \beta > 2 + 4/\rho$, *then for every* $k \in \mathbf{Z}_+^d$ *with* $L_k \neq \varnothing$, *we have*

$$\sum_{i \in L_k} |v_i| \leq C[N_k]^{1/2-\alpha_0} \qquad a.s. \quad and \quad \sum_{i \in L_k} |w_i| \leq C[N_k]^{1/2-\alpha_0} \qquad a.s.$$

PROOF. For the first inequality, we follow the proof of Lemma 8 of Berkes and Morrow (1981). Note that $I_i = \bigcup_{s=1}^d I_i(s)$, where $I_i(s)$ are disjoint rectangles with $|I_i(s)| \leq Ci_s^\beta \prod_{s' \neq s} i_{s'}^\alpha$. Hence, $v_i = \sum_{s=1}^d v_i(s)$ with $v_i(s) := \sum_{j \in I_i(s)} X_j$.

By Chebyshev's inequality and (4),

$$P(|v_i(s)| \geq [i]^{\alpha/2-1}) \leq C[i]^{-(\alpha-2)}|I_i(s)| \leq Ci_s^{-(\alpha-\beta-2)}\prod_{s' \neq s} i_{s'}^2$$
$$\leq i_s^{-(\alpha-\beta-2-2/\rho)} \leq C[i]^{-(\alpha-\beta-2-2/\rho)\rho/2}$$

for every $i \in L_k$. (As in the proof of the above-mentioned lemma, we used the fact that $i \in L_k$ implies that $i_s \geq C \prod_{s' \neq s} i_{s'}^\rho$ and, consequently, $i_s \geq$ $C[i]^{\rho/2}$.) Since $(\alpha - \beta - 2 - 2/\rho)\rho/2 > 1$, the result follows by the Borel–Cantelli lemma. A similar argument applies to $w_i$, since $E(w_i^2) = \tau_i^2 \leq C|I_i| = C\sum_{s=1}^{d}|I_i(s)|$. □

**4. The approximation theorem.** In this section we will verify that the third sum in the decomposition (6) of $S(R_k)$ can be approximated by $W(R_k)$, where $W$ is a $d$-parameter Wiener process with variance $\sigma^2$. Some preliminary lemmas are needed.

The next result follows exactly as Theorem 2.1 of Yu (1996), using Lemma 3.2.

LEMMA 4.1. *If* (C1) *and* (C2) *hold and* $2r_0r/(2+r) < \alpha/\beta < 2(1+r)/(2+r)$ *with* $r_0 := \max\{1, (r+\delta)^{-1}\}$, *then for any* $0 < \theta < 1/2$ *and all* $i \neq j$,

$$E(\eta_i \eta_j) \leq C\{([i][j])^{-\alpha/2} E(u_i u_j)\}^{\theta/(1+\theta)}.$$

The next lemma gives a generalization of relationship (3.11) of Yu (1996) to the multi-parameter case.

LEMMA 4.2. *If* (C2) *holds, then*

$$E(u_i u_j) \leq C e^{-\lambda M_{i,j}^{\beta}},$$

*where* $M_{i,j} := \max_{s\,:\,i_s \neq j_s}(M_s(i,j) - 1)$ *and* $M_s(i,j) := \max(i_s, j_s), s = 1, \ldots, d$.

PROOF. Let $d := \min_{k \in H_i} d(k, H_j)$ be the distance between $H_i$ and $H_j$, where $d(k, H_j) := \min_{k' \in H_j} \|k - k'\|$. Then $d_k := d(k, H_j) - d \geq 0 \; \forall k \in H_i$,

$$E(u_i u_j) = \sum_{k \in H_i}\sum_{k' \in H_j} E(X_k X_{k'}) \leq \sum_{k \in H_i} u(d + d_k) \leq Ce^{-\lambda d}\sum_{k \in H_i} e^{-\lambda d_k} \leq Ce^{-\lambda d}$$

and $d = \max_{s=1,\ldots,d} \min_{k \in H_i, k' \in H_j}|k_s - k'_s| = \max_{s\,:\,i_s \neq j_s}\{m_s^\beta + \sum_{l=m_s+1}^{M_s-1}(l^\alpha + l^\beta)\} \geq M_{i,j}^\beta$, where $m_s = m_s(i,j) := \min(i_s, j_s)$ and $M_s = M_s(i,j)$. □

In order to prove our approximation theorem, we need to be able to "count" properly the indices in $L$, that is, to define a bijection $\psi: \mathbf{Z}_+ \to L$ satisfying certain properties. This will be given by the following lemma, whose proof can be found in the Appendix.

LEMMA 4.3. *There exists a bijection* $\psi: \mathbf{Z}_+ \to L$ *such that*

(8) $\quad l < m \implies \exists s^* = s^*(l,m) \text{ such that } \psi(l)_{s^*} \leq \psi(m)_{s^*}$

(9) $\quad \exists m_0 \in \mathbf{Z}_+ \text{ such that } m \leq C[\psi(m)]^{\gamma_0} \quad \forall m \geq m_0$

*for any* $\gamma_0 > (1 + 1/\rho)(1 - 1/d)$.

We are now able to prove the desired approximation theorem.

THEOREM 4.4. *Suppose that* (C1) *and* (C2) *hold,* $\alpha > 3(1 + 1/\rho)(1 - 1/d)$, $\beta > (2/\rho)(1 + 1/\rho)(1 - 1/d)$ *and* $2r_0 r/(2+r) < \alpha/\beta < 2(1+r)/(2+r)$ *with* $r_0 := \max\{1, (r+\delta)^{-1}\}$. *Then without changing its distribution, we can redefine the random field* $(X_j)_{j \in \mathbf{Z}_+^d}$ *on a rich enough probability space together with a* $d$*-parameter Wiener process* $W = (W_t; t \in [0, \infty)^d)$ *with variance* $\sigma^2$, *such that for every* $k \in \mathbf{Z}_+^d$ *with* $L_k \neq \varnothing$,

$$\sum_{i \in L_k} \sigma\sqrt{|B_i|}\left|\eta_i - \frac{W(B_i)}{\sigma\sqrt{|B_i|}}\right| \leq C[N_k]^{1/2 - \alpha_0} \qquad a.s.,$$

*where* $\alpha_0 := 1/\{2(1+\alpha)\}$.

PROOF. Let $0 < \theta < 1/2$ be such that $\alpha\{(1+1/\rho)(1-1/d)\}^{-1} > 1 + 1/\theta$ and choose $\gamma_0$ such that $(1 + 1/\rho)(1 - 1/d) < \gamma_0 < \min\{\alpha\theta/(1+\theta), \beta\rho/2\}$. Let $\psi: \mathbf{Z}_+ \to L$ be the bijection given by Lemma 4.3.

We will apply Theorem 5 of Berkes and Philipp (1979) to the sequence $Y_m := \eta_{\psi(m)}, m \in \mathbf{Z}_+$ of random variables and the probability distributions $G_m := N(0, 1), m \in \mathbf{Z}_+$ and we will prove that for each $m \in \mathbf{Z}_+, m \geq 2$ there exists some $\rho_m > 0$ such that

$$(10) \qquad \left| E \exp\left\{i \sum_{l=1}^m t_l Y_l\right\} - E \exp\left\{i \sum_{l=1}^{m-1} t_l Y_l\right\} E \exp\{it_m Y_m\}\right| \leq \rho_m$$

for all $t_1, \ldots, t_m \in \mathbf{R}$ with $\sum_{l=1}^m t_l^2 \leq U_m^2$, where $U_m > 10^4 m^2$.

Then, by the above-mentioned theorem, without changing its distribution, we can redefine the sequence $(Y_m)_{m \in \mathbf{Z}_+}$ on a rich enough probability space together with a sequence $(Z_m)_{m \in \mathbf{Z}_+}$ of independent $N(0, 1)$-random variables such that

$$P(|Y_m - Z_m| \geq \alpha_m) \leq \alpha_m \qquad \forall m \in \mathbf{Z}_+,$$

where $\alpha_m \leq C\{U_m^{-1/4} \log U_m + \exp(-3U_m^{1/2}/16)m^{1/2}U_m^{1/4} + \rho_m^{1/2}U_m^{m+1/4}\}$. We will prove next that

$$(11) \qquad \qquad \alpha_m \leq Cm^{-2} \qquad \text{for } m \text{ large}.$$

Then, by the Borel–Cantelli lemma, $|Y_m - Z_m| \leq C\alpha_m$, $\forall m \in \mathbf{Z}_+$ a.s. Using a straightforward $d$-parameter generalization of Lemma 4 of Csörgő and Révész (1975), without changing its distribution, we can redefine the sequence $(Z_m)_{m \in \mathbf{Z}_+}$ on a richer probability space together with a $d$-parameter Wiener process with variance $\sigma^2$ such that $Z_m = W(B_{\psi(m)})/(\sigma\sqrt{|B_{\psi(m)}|})$, $\forall m \in \mathbf{Z}_+$. Hence,

$$\left|\eta_i - \frac{W(B_i)}{\sigma\sqrt{|B_i|}}\right| \leq C\alpha_{\psi^{-1}(i)} \qquad \forall i \in L \text{ a.s.}$$

and because $|B_i| \leq |B_k| \leq C[k]^\alpha$, $\forall\, i \in L_k$ and $\sum_{l \in \mathbf{Z}_+} \alpha_l < \infty$, we have

$$\sum_{i \in L_k} \sigma\sqrt{|B_i|}\left|\eta_i - \frac{W(B_i)}{\sigma\sqrt{|B_i|}}\right| \leq C[k]^{\alpha/2} \sum_{i \in L_k} \alpha_{\psi^{-1}(i)} \leq C[k]^{\alpha/2} \leq C[N_k]^{1/2-\alpha_0}.$$

We proceed next to the verification of (10) and (11). By Lemmas 4.1 and 4.2, we have

$$E(Y_l Y_m) \leq C\{([\psi(l)][\psi(m)])^{-\alpha/2} E(u_{\psi(l)} u_{\psi(m)})\}^{\theta/(1+\theta)}$$

$$\leq C([\psi(l)][\psi(m)])^{-\alpha\theta/(2+2\theta)} e^{-\lambda\theta M^\beta_{\psi(l),\psi(m)}/(1+\theta)}$$

$$\leq C([\psi(l)][\psi(m)])^{-\alpha\theta/(2+2\theta)} e^{-\lambda\theta[\psi(m)]^{\beta\rho/2}/(1+\theta)}.$$

[For the last inequality above we used (8) to obtain an $s^* = s^*(l,m)$, for which $M_{s^*}(\psi(l),\psi(m)) = \psi(m)_{s^*}$; since $\psi(m) \in L$, we have $M_{\psi(l),\psi(m)} \geq \psi(m)_{s^*} - 1 \geq C[\psi(m)]^{\rho/2}$.] By Lemma 2.2 of Dabrowski and Dehling (1988), the left-hand side of (10) is smaller than $2\sum_{l=1}^{m-1} |t_l t_m| E(Y_l Y_m)$, which is, in turn, smaller than

$$Ce^{-\lambda\theta[\psi(m)]^{\beta\rho/2}/(1+\theta)} \sum_{l=1}^{m-1} 2|t_l t_m|([\psi(l)][\psi(m)])^{-\alpha\theta/(2+2\theta)}$$

$$\leq Ce^{-\lambda\theta[\psi(m)]^{\beta\rho/2}/(1+\theta)} \left\{\sum_{l=1}^{m-1} t_l^2 [\psi(l)]^{-\alpha\theta/(1+\theta)}\right.$$

$$\left. + (m-1)t_m^2 [\psi(m)]^{-\alpha\theta/(1+\theta)}\right\}$$

$$\leq Ce^{-\lambda\theta[\psi(m)]^{\beta\rho/2}/(1+\theta)} \sum_{l=1}^{m} t_l^2 \leq Ce^{-\lambda\theta[\psi(m)]^{\beta\rho/2}/(1+\theta)} U_m^2 := \rho_m$$

for $m$ large enough. (In the second inequality above, we used the fact that $m \leq C[\psi(m)]^{\alpha\theta/(1+\theta)}$, which follows from Lemma 4.3 by our choice of $\gamma_0$.)

Finally, relationship (11) follows if we take $U_m := m^q$ with $q > 8$. Clearly, $U_m^{-1/4} \log U_m \leq m^{-2}$ and $\exp(-3U_m^{1/2}/16)m^{1/2}U_m^{1/4} \leq \exp(-2U_m^{1/2}/16) \leq m^{-2}$ for $m$ large enough. We have

$$\rho_m^{1/2} U_m^{m+1/4} = e^{-\lambda\theta[\psi(m)]^{\beta\rho/2}/(2+2\theta)} m^{q(m+5/4)} \leq m^{-2}$$

since $\{2 + q(m+5/4)\}\log m \leq Cm^{1+\varepsilon} \leq C[\psi(m)]^{(1+\varepsilon)\gamma_0} \leq C[\psi(m)]^{\beta\rho/2}$, for $m$ large enough. This concludes the proof of the theorem. $\square$

REMARK. A similar argument can be used to give a simplified proof for Theorem 2.5 of Yu (1996) (in the case $d = 1$). More precisely, one can check directly the condition of Theorem 5 of Berkes and Philipp (1979) for the sequence $(\eta_k)_{k\geq 1}$ of random variables and the probability distributions $G_k = N(0,1), k \geq 1$ (as we did above). We obtain in this manner a sequence $(Z_k)_{k\geq 1}$ of independent $N(0,1)$-random variables with $P(|\eta_k - Z_k| \geq \alpha_k) \leq \alpha_k$ and $\alpha_k \leq Ck^{-2}$. Without changing its distribution, we can redefine the sequence $(Z_k)_{k\geq 1}$ on a richer probability space together with a standard Brownian motion $W = \{W_t; t \in [0, \infty)\}$ such that $Z_k = W(\hat{H}_k)/\sqrt{\lambda_k^2 + \tau_k^2}$, where $\hat{H}_k := (V_{k-1}, V_k]$ and $V_k := \sum_{i=1}^k (\lambda_i^2 + \tau_i^2)$. Since $\lambda_i^2 + \tau_i^2 \leq Ci^\alpha \leq Ck^\alpha$ for $i \leq k$ and $\sum_{i\geq 1} \alpha_i < \infty$, this gives immediately the desired approximation

$$\sum_{i=1}^k \sqrt{\lambda_i^2 + \tau_i^2} \left| \eta_i - \frac{W(\hat{H}_i)}{\sqrt{\lambda_i^2 + \tau_i^2}} \right| \leq Ck^{\alpha/2} \sum_{i=1}^k \alpha_i \leq CN_k^{1/2-\alpha_0} \quad \text{a.s.}$$

**5. The remaining terms.** In this section we show that the terms $S((0, N_k] \setminus R_k), W((0, N_k] \setminus R_k), S_N - S_{N_k}, W_N - W_{N_k}$ can be made sufficiently small if $N \in G_\tau$.

Note that $(0, N_k] \setminus R_k = \bigcup_{s=1}^d (0, N_k^{(s)}]$. If we let $D_s(N) := \max_{n \leq N^{(s)}} |S_n|$ and $\hat{D}_s(N) := \max_{n \leq N^{(s)}} |W_n|$, for each $s = 1, \ldots, d$ and $N \in H$, then

$$S((0, N_k] \setminus R_k) \leq \sum_{s=1}^d 2^{d-s} D_s(N_k), \qquad W((0, N_k] \setminus R_k) \leq \sum_{s=1}^d 2^{d-s} \hat{D}_s(N_k).$$

On the other hand, $(0, N] \setminus (0, N_k] = \bigcup_J I_k^{(J)}$, where $I_k^{(J)} := \prod_{s \in J} (n_{k_s}, N_s] \times \prod_{s \in J^c} (0, n_{k_s}]$ and the union is taken over all nonempty subsets $J$ of $\{1, \ldots, d\}$. Let $M_k^{(J)} := \max |S(I_k^{(J)})|$ and $\hat{M}_k^{(J)} := \sup |W(I_k^{(J)})|$, where the maximum and the supremum are taken over all $N$ with $n_{k_s} < N_s \leq n_{k_s+1}$, $\forall s \in J$. We have

$$\max_{N_k < N \leq N_{k+1}} |S_N - S_{N_k}| \leq \sum_J M_k^{(J)}, \qquad \sup_{N_k < N \leq N_{k+1}} |W_N - W_{N_k}| \leq \sum_J \hat{M}_k^{(J)}.$$

We note in passing that the arguments that are valid for the terms depending on the original random field $(X_j)_{j \in \mathbf{Z}_+^d}$ can be applied to the terms depending on the Wiener process $W$, since $W(V) = \sum_{j \in V} \hat{X}_j$, $\forall V \in \mathcal{A}$, where $\hat{X}_j := W((j-1, j])$ are independent $N(0, \sigma^2)$-random variables. Clearly, $(\hat{X}_j)_{j \in \mathbf{Z}_+^d}$ is a weakly stationary associated random field with zero mean and covariance coefficient $\hat{u}(n) = 0$, $\forall n \geq 1$.

LEMMA 5.1. (a) *Suppose that* (C1) *and* (C2′) *hold with* $\nu \geq d\nu_0$ *and* $\nu_0 := r(2 + r + \delta)/(2\delta) < (d-2)^{-1}$ *if* $d \geq 3$. *Then there exists* $x_0$ *such that* $\forall V \in \mathcal{A}$, $\forall x \geq x_0$,

$$P(M(V) \geq x|V|^{1/2}) \leq Cx^{-(2+r)},$$

INVARIANCE PRINCIPLE FOR RANDOM FIELDS                    15

where $M(V) := \max\{|S(Q)|; Q \subseteq V, Q \in \mathcal{A}\}$.

(b) *If* (C1) *and* (C2) *hold, then there exists* $\gamma > 0$ *such that* $\forall V \in \mathcal{A}$,

$$P(\tilde{M}(V) \geq |V|^{1/2}(\log |V|)^{d+1}) \leq C|V|^{-\gamma},$$

where $\tilde{M}((a,b]) := \max\{|S(Q)|; Q = (a,c], a < c \leq b\}$.

PROOF. (a) Using Lemma 1 of Bulinski and Keane (1996), the Markov inequality and Lemma 3.4, we have $P(M(V) \geq x|V|^{1/2}) \leq 2P(|S(V)| \geq x|V|^{1/2}/2) \leq Cx^{-(2+r)}|V|^{-(1+r/2)}E|S(V)|^{2+r} \leq Cx^{-(2+r)}$.

(b) This follows exactly as the second inequality of Lemma 7 of Berkes and Morrow (1981), using the moment inequality given by Lemma 3.4 and the rate of convergence in the CLT given by Lemma 3.1. This rate is sharper than the rate of Lemma 5 of Berkes and Morrow (1981). To see this, we use (4) and we note that $\sup_{x \in \mathbf{R}} |F_V(x) - \Phi(x)|$ is either smaller than $C|V|^{-\{s/2-1-\varepsilon d(s-1)\}}$ if $s \leq 3$, or smaller than $C|V|^{-(1/2-\varepsilon d)}$ if $s > 3$; in both cases a suitable choice of $\varepsilon > 0$ gives us the rate $C|V|^{-t}$ for some $t \in (0,1)$. We also note that the requirement $|V| \in G_\tau$ is not needed. □

The next result follows exactly as Lemma 6 of Berkes and Morrow (1981), using Lemma 5.1(a).

LEMMA 5.2. *If* $\alpha > 16/(3\tau) - 1$, *then*

$$\max_{s=1,\ldots,d} D_s(N_k) \leq C[N_k]^{1/2-\varepsilon} \qquad a.s.,$$

$$\max_{s=1,\ldots,d} \hat{D}_s(N_k) \leq C[N_k]^{1/2-\varepsilon} \qquad a.s.$$

*for every* $N_k \in G_\tau$ *and* $0 < \varepsilon < \tau/32$.

The following result follows exactly as Lemma 9 of Berkes and Morrow (1981), using Lemma 5.1(b).

LEMMA 5.3. *Let* $\gamma$ *be the constant given by Lemma* 5.1(b). *If* $\alpha > 2/\gamma$, *then*

$$\max_J M_k^{(J)} \leq C[N_k]^{1/2-\varepsilon} \qquad a.s.,$$

$$\max_J \hat{M}_k^{(J)} \leq C[N_k]^{1/2-\varepsilon} \qquad a.s.$$

*for every* $N_k \in G_\rho$ *and* $0 < \varepsilon < \rho/(8\alpha)$.

## APPENDIX

PROOF OF LEMMA 3.2. Using Lemma 3.1 for $V = H_k$ and relationship (5), we obtain that $\sup_{x \in \mathbf{R}} |\tilde{F}_k(x) - \Phi(x)|$ is either smaller than $C[k]^{-\{\alpha s/2 - \alpha - \varepsilon d(s-1)\}}$ if $s \leq 3$, or smaller than $C[k]^{-(\alpha/2 - \varepsilon d)}$ if $s > 3$. If $\alpha/\beta > 2r_0 r/(2+r)$, then a suitable choice of $\varepsilon > 0$ allows us to conclude that $|\tilde{F}_k(x) - \Phi(x)| \leq C[k]^{-r\beta/(2+r)}$, $\forall x \in \mathbf{R}$. The first inequality follows by a change of variables.

For the second inequality we use a technique similar to that used to prove relationship (3.3) of Yu (1996). Let $\varphi_k(t) := E[\exp(it\xi_k)], \tilde{\varphi}_k(t) := E[\exp(itu_k/\lambda_k)]$ and $\varphi(t) = \exp(-t^2/2)$. Since $(\lambda_k^2 + \tau_k^2)/\lambda_k^2 \leq C$, we have, for any $T > 0$,

$$
\begin{aligned}
|f_k(x) - f(x)| &\leq \frac{1}{2\pi} \int_{-\infty}^{\infty} |\varphi_k(t) - \varphi(t)| \, dt \\
&\leq \frac{C}{2\pi} \int_{-\infty}^{\infty} |\tilde{\varphi}_k(t) - \varphi(t)| \exp\left\{-\frac{\tau_k^2 t^2}{2\lambda_k^2}\right\} ds \\
&\leq \frac{C}{2\pi} \cdot 2T[k]^{-r\beta/(2+r)} + \frac{C}{\pi} \int_{|t| \geq T} \exp\left\{-\frac{t^2 \tau_k^2}{2\lambda_k^2}\right\} dt \\
&\leq C \cdot T[k]^{-r\beta/(2+r)} + \frac{C}{T} \cdot \frac{\lambda_k^2}{\tau_k^2} \exp\left\{-\frac{\tau_k^2 T^2}{2\lambda_k^2}\right\}.
\end{aligned}
$$

Since $\lambda_k^2/\tau_k^2 \leq C[k]^{\alpha-\beta}$, the conclusion follows by choosing $T = C[k]^q$ with $\alpha - \beta < q < r\beta/(2+r)$. Such a choice is possible if $\alpha/\beta < 2(1+r)/(2+r)$. □

PROOF OF LEMMA 3.7. First we claim that it is enough to prove (7) for "squares," that is, for rectangles $V = (m, n] \in \mathcal{A}$ for which $n_s - m_s = l$, $\forall s = 1, \ldots, d$. To see this, we note that each rectangle $V$ can be written as a finite union of disjoint squares: $V = \bigcup_{i=1}^{p} V_i$. By the association property $\sigma^2(V) \geq \sum_{i=1}^{p} \sigma^2(V_i)$ and

$$\sigma^2 - \frac{\sigma^2(V)}{|V|} \leq \frac{1}{|V|} \sum_{i=1}^{p} |V_i| \left(\sigma^2 - \frac{\sigma^2(V_i)}{|V_i|}\right) \leq \frac{1}{|V|} \sum_{i=1}^{p} C|V_i|^{1-\delta_0} \leq C|V|^{-\delta_0}$$

because $0 < \delta_0 < 1$. Let us now prove relationship (7) for a square $V = (m, n]$ with $n_s - m_s = l$, $\forall s = 1, \ldots, d$. Note that $|V| = l^d$. By stationarity,

$$
\begin{aligned}
\sigma^2(V) &= |V| \cdot r(0) + \sum_{-(n-m-1) \leq i \leq n-m-1, i \neq 0} \prod_{s=1}^{d} (l - |i_s|) \cdot r(i) \\
&= |V| \cdot \sum_{\|i\| \leq l-1} r(i) - \sum_{\varnothing \neq K \subseteq \{1,\ldots,d\}} (-1)^{|K|-1} \sum_{\|i\| \leq l-1, i \neq 0} c(K, i) \cdot r(i),
\end{aligned}
$$

where $c(K,i) := l^{|K^c|} \cdot \prod_{s \in K} |i_s|$. Since $\sigma^2 - \sum_{\|i\| \leq l-1} r(i) = \sum_{\|i\| \geq l} r(i) = u(l)$ and $c(K,i) \leq |V|$ if $\|i\| \leq l-1$, we have

$$\sigma^2 - \frac{\sigma^2(V)}{|V|} \leq u(l) + \sum_{\varnothing \neq K \subseteq \{1,\ldots,d\}, |K| \text{ odd}} \frac{1}{|V|} \sum_{\|i\| \leq l-1, i \neq 0} c(K,i) \cdot r(i)$$

$$\leq C|V|^{-\nu/d} + \sum_{\varnothing \neq K \subseteq \{1,\ldots,d\}, |K| \text{ odd}} \sum_{\|i\| \leq l-1, i \neq 0} r(i)$$

$$\leq C|V|^{-\nu/d} + C|V|^{-\nu/d+1}.$$

We used the fact that $u(l) \leq Cl^{-\nu} = C|V|^{-\nu/d}$ and $r(i) \leq u(\|i\|) \leq u([i]^{1/d}) \leq C[i]^{-\nu/d}$ for any $i \in \mathbf{Z}^d$, where $[i] = \prod_{s: i_s \neq 0} |i_s|$. $\square$

PROOF OF LEMMA 4.3. The idea of the proof is based on the following simple observation in the case $d = 2$. For each $m \in \mathbf{Z}_+, m \geq 2$ with $(m,m) \in L$, there exists a $k_1^*(m) \geq m$ such that $(k_1, m), (m, k_1) \in L$ for every $m \leq k_1 \leq k_1^*(m)$. Therefore, to each vertex $(m,m) \in L$, one can associate an "L-shaped" region $L(m)$ consisting of $2\{k_1^*(m) - m\} + 1$ points in $L$. In view of the desired property (8), we will count consecutively the indices in $L(2), L(3)$, and so on. To verify property (9), we note that $k \in L(m)$ implies $[k] \geq m^2$.

We begin now the proof for arbitrary $d \geq 2$. Let $m \in \mathbf{Z}_+, m \geq 2$ be such that $(m, \ldots, m) \in L$ and $k = (k_1, \ldots, k_{d-1}, m) \in L$ be such that $k_s > m, \forall s < d$. This implies that all the vertices of $B_k$ are in $G_\rho$, and, in particular, $n_m \geq n_{k_s}^\rho, \forall s < d$. Since $m$ is fixed, this cannot happen for infinitely many $k_s$'s. It follows that for each $s = 1, \ldots, d-1$, there exists a $k_s^*(m) \geq m$ such that $k_s \leq k_s^*(m)$. We note that $k_s^*(m) \leq Cm^{1/\rho}$, if $m$ is large enough. This argument shows us that we have a maximum number of $k^*(m) := \prod_{s=1}^{d-1}\{k_s^*(m) - m\}$ points of the form $(k_1, \ldots, k_{d-1}, m)$ in $L$, with $k_s > m, \forall s < d$.

By symmetry, we can repeat this argument for each of the axes. We let $L_s(m) := \{k = (k_1, \ldots, k_{s-1}, m, k_s, \ldots, k_{d-1}); m < k_{s'} \leq k_{s'}^*(m), \forall s' < d\}$ for every $s = 1, \ldots, d$. The "L-shaped" region corresponding to the index $m$ is

$$L(m) := \bigcup_{s=1}^{d} L_s(m) \cup \{(m, \ldots, m)\}.$$

Note that $|L(m)| = dk^*(m) + 1$ and that $k \in L(m)$ implies $[k] \geq m^d$. Clearly, $L \subseteq \bigcup_m L(m)$ [note that in the case $d = 2$, we actually have $L = \bigcup_m L(m)$].

Next we count consecutively the indices in $L(2), L(3)$, and so on, that is, we define a bijection $\varphi: \mathbf{Z}_+ \to \bigcup_m L(m)$ such that $\forall z \in \mathbf{Z}_+$,

$$\sum_{l=2}^{m-1} |L(l)| < z \leq \sum_{l=2}^{m} |L(l)| \implies \varphi(z) \in L(m).$$

The bijection $\varphi$ clearly satisfies condition (8). To verify (9), we note that

$$z \leq d \sum_{l=2}^{m} \prod_{s=1}^{d-1} (k_s^*(l) - l) + m \leq d \prod_{s=1}^{d-1} \sum_{l=2}^{m} (k_s^*(l) - l) + m$$

$$\leq d \prod_{s=1}^{d-1} \left( \sum_{l=2}^{m} k_s^*(l) \right) - d \left( \sum_{l=2}^{m} l \right)^{d-1} + m \leq d \prod_{s=1}^{d-1} \left( \sum_{l=2}^{m} k_s^*(l) \right)$$

$$\leq C m^{(1+1/\rho)(d-1)} \leq C m^{d \gamma_0} \leq C [\varphi(z)]^{\gamma_0}$$

for $m$ large enough and $\gamma_0 > (1+1/\rho)(1-1/d)$ arbitrary. Finally, define the bijection $\psi : \mathbf{Z}_+ \to L$ such that $\psi^{-1}(k) \leq \varphi^{-1}(k)$, $\forall k \in L$. The result follows since if $z_1, z_2 \in \mathbf{Z}_+$ are such that $\psi(z_1) = \varphi(z_2)$, then $z_1 \leq z_2$. $\square$

Department of Mathematics and Statistics
University of Ottawa
585 King Edward Avenue
Ottawa, Ontario
Canada K1N 6N5
e-mail: rbala348@science.uottawa.ca
url: http://aix1.uottawa.ca/~rbalan